# Least Square Estimation-Based SDP Cuts for SOCP Relaxation of AC OPF

Zhixin Miao, *Senior Member, IEEE,* Lingling Fan, *Senior Member, IEEE*, Hoseein Ghassempour, and Bo Zeng, *Member, IEEE*

*Abstract*—This paper presents a method that generates affine inequalities to strengthen the second-order conic programming (SOCP) relaxation of an alternating current optimal power flow (AC OPF) problem. The affine inequalities serve as cuts to get rid of points outside of the feasible region of AC OPF with semi-definite programming (SDP) relaxation. Hence, the affine inequalities are names as SDP cuts. While AC OPF with SDP relaxation has a high computational complexity, AC OPF with SOCP has a much lower computational complexity. Recent research has found that the feasible region of SDP relaxation is contained inside the feasible region of the SOCP relaxation. Therefore, the integration of SDP cuts into SOCP relaxation provides better scalability compared to the SDP relaxation and a tighter gap compared to the SOCP relaxation. The SDP cuts are generated by solving least square estimation (LSE) problems at cycle basis and further exploring the geometric characteristic of LSE. General feasibility cuts generating method is also employed for analysis. We found that the SDP cuts generated by LSE method are indeed feasibility cuts. The SDP cuts effectively reduce the search space. Case studies of systems with several buses to hundreds of buses have demonstrated the method is very effective in reducing the gaps.

## I. INTRODUCTION

SDP relaxation and SOCP relaxation have been applied in AC OPF in the literature [1], [2]. While AC OPF with SDP relaxation has scalability issues, AC OPF with SOCP has demonstrated excellent scalability. Recent research has found that the feasible region of SDP relaxation is contained inside the feasible region of the SOCP relaxation [3], [4]. The standard SOCP relaxation does not take care of the meshed network cycle constraints (sum of the voltage angle differences across a cycle is zero). Most recently, [3] presents three approaches to tackle the cycle issue in AC OPF with SOCP relaxation. One of the approaches is to generate SDP cuts represented by linear inequalities. This approach is claimed to produce tighter gaps compared to other approaches.

Most recently, [5] proposed to use nonlinear inequalities to generate SDP cuts for SOCP AC OPF. The inequalities are based on the matrix determinant. For large networks, nonlinear inequalities are expensive for computation. Compared to nonlinear inequalities, linear or affine inequalities will have more computational advantages for large-size networks.

In this paper, we will examine the method presented in [3], show its deficiency, and present a different method to generate

Z. Miao, L. Fan, and H Ghassempour are with Department of Electrical Engineering at University of South Florida, Tampa, FL 33620 (Email: linglingfan@usf.edu). B. Zeng is with University of Pittsburgh.

SDP cuts. In our paper, the SDP cuts represented by affine inequalities will be generated by least square estimation (LSE). We then compare the LSE-based method with the general feasibility cutting plane method and demonstrate the similar nature of two methods. Test cases of several buses to hundreds of buses are examined. With a few iterations, the approach will generate a lower bound with much tighter gap compared to that generated by SOCP relaxation.

The rest of the paper is organized as follows. Section II presents AC OPF relaxations. Both SOCP and SDP relaxations are presented. The SDP separation method presented in [3] is examined in details. Section III presents the proposed LSE method and feasibility cut method. Section IV presents case study results and Section V concludes the paper.

## II. AC OPF RELAXATIONS

AC OPF is formulated as an optimization problem with the objective function to minimize the cost of generation or power loss, equality constraints representing the relationship of bus power injection versus bus voltage magnitudes (notated by a vector $V \in \mathbb{R}^n$, where $n$ is the number of buses in the system) and phase angles ($\theta \in \mathbb{R}^n$), and inequality constraints representing voltage limits, line flow limits, generation limits, etc [6]. The decision variables of AC OPF are voltage magnitudes $V$, phase angles $\theta$, and generators' real and reactive power outputs, notated as $P_g \in \mathbb{R}^{|\mathcal{G}|}$ and $Q_g \in \mathbb{R}^{|\mathcal{G}|}$, where $\mathcal{G}$ is the set of generators and $|.|$ notates the cardinality of a set.

AC OPF is a nonconvex optimization problem. This can be shown by the power injection equality constraints as follows. Given the system admittance matrix $Y = G + jB$, the power injection at every node can be expressed by $V$ and $\theta$.

$$P_i^g - P_i^d = \sum_{j=1}^n V_i V_j (G_{ij} \cos(\theta_i - \theta_j) + B_{ij} \sin(\theta_i - \theta_j))$$

$$Q_i^g - Q_i^d = \sum_{j=1}^n V_i V_j (G_{ij} \sin(\theta_i - \theta_j) - B_{ij} \cos(\theta_i - \theta_j))$$

where superscript $g$ notates generator's output and $d$ notates load consumption.

Note that the equality constraints of power injections are non-convex in terms of $V$ and $\theta$. Relaxations have been developed in the literature to have a convex feasible region. These methods deal with new sets of decision variables to replace $V$ and $\theta$.

## A. SOCP relaxation

For AC OPF, a new set of variables $c_{ij}$ and $s_{ij}$ are used to replace the voltage phasors $V_i \angle \theta_i, i \in \mathcal{B}$, where $\mathcal{B}$ is the set of buses in a power network.

$$c_{ii} = V_i^2, \ c_{ij} = V_i V_j \cos(\theta_i - \theta_j)$$
$$s_{ii} = 0, \ s_{ij} = -V_i V_j \sin(\theta_i - \theta_j) \quad (1)$$

where $c_{ij} = c_{ji}$ and $s_{ij} = -s_{ji}$.

It is easy to find the following relationship:

$$c_{ij}^2 + s_{ij}^2 = V_i^2 V_j^2 = c_{ii} c_{jj}. \quad (2)$$

The AC OPF problem's power injection constraints are linear in terms of $c_{ij}$ and $s_{ij}$.

$$P_i^g - P_i^d = \sum_{j=1}^n (G_{ij} c_{ij} - B_{ij} s_{ij}), \quad (3)$$

$$Q_i^g - Q_i^d = \sum_{j=1}^n (-G_{ij} s_{ij} - B_{ij} c_{ij}) \quad (4)$$

where $n$ is the total number of buses.

In addition, (2) will be relaxed as a second-order cone:

$$c_{ij}^2 + s_{ij}^2 \leq c_{ii} c_{jj}. \quad (5)$$

The above relaxation is named as SOCP relaxation. Research work has been conducted in this area and sufficient conditions for SOCP relaxation being exact are also found [7], [8]. In spanning tree power networks, under a mild condition (e.g., voltage upper bounds not binding), SOCP relaxations are exact [8]. For meshed network, SOCP relaxation is not exact since the following constraint is not considered in SOCP relaxation:

$$\tan(\theta_i - \theta_j) = -\frac{s_{ij}}{c_{ij}} \quad (6)$$

For a cycle $C$, the angle constraint is as follows

$$\sum_{(i,j) \in C} \theta_{ij} = 2\pi k, \text{for some } k \in \mathbb{Z}. \quad (7)$$

Without imposing constraint (6), it is possible to end up with $c_{ij}, s_{ij}$ that violate the cycle constraint, or infeasible solution to the original AC OPF problem.

To solve this issue, in [9] Jabr proposed linear approximation for (6). This method requires iteration. In addition, feasible points could be lost due to the imposed linear constraint. Kocuk *et al* proposed three methods in [3], including McCormick based linear programming relaxation and separation, SDP separation, and arctangent envelopes to deal with the cycle constraints. In [3], the authors claim that the SDP separation works the best in term of providing tight gaps. The SDP cuts in [3] are represented by linear inequalities.

Most recently, [5] proposed to use nonlinear inequalities to generate SDP cuts for SOCP AC OPF. The inequalities are based on the matrix determinant. For large networks, nonlinear inequality may be expensive in computation.

Compared to nonlinear inequalities, linear or affine inequalities will have more computational advantages for large-size networks.

In this article, we will examine the SDP separation in [3] and develop SDP cuts using least square estimation (LSE). As a first task, the relationship of decision variables in SOCP relaxation $c_{ij}, s_{ij}$ and those in SDP relaxation is examined. Then we design the algorithm to generate SDP cuts.

## B. SDP relaxation

In SDP relaxation, rectangular expressions are used to represent the voltage phasors.

$$V_i = V_i \angle \theta_i = \underbrace{V_i \cos \theta_i}_{e_i} + j \underbrace{V_i \sin \theta_i}_{f_i} = e_i + j f_i \quad (8)$$

A matrix $W$ is defined as follows

$$W = \begin{bmatrix} e \\ f \end{bmatrix} \begin{bmatrix} e^T & f^T \end{bmatrix} \quad (9)$$

where $f = (f_1, f_2, \cdots, f_n)^T$, $e = (e_1, e_2, \cdots, e_n)^T$.

It's obvious to find the following characteristics:

$$W = W^T \text{ and } W \succeq 0 \quad (10)$$

$W \succeq 0$ means that this matrix is semi-definite.

The power injection constraints will be shown to be linear with the elements of $W$. Define

$$i' = i + |\mathcal{B}|, \quad j' = j + |\mathcal{B}|, \quad (11)$$

where $|.|$ notates the cardinality of a set.

$$P_i^g - P_i^d = \sum_{j=1}^n V_i V_j (G_{ij} \cos(\theta_i - \theta_j) + B_{ij} \sin(\theta_i - \theta_j))$$
$$= \sum_{j=1}^n (G_{ij}(e_i e_j + f_i f_j) + B_{ij}(f_i e_j - e_i f_j))$$
$$= \sum_{j=1}^n G_{ij}(W_{ij} + W_{i'j'}) + B_{ij}(W_{ji'} - W_{ij'})$$

$$Q_i^g - Q_i^d = \sum_{j=1}^n V_i V_j (G_{ij} \sin(\theta_i - \theta_j) - B_{ij} \cos(\theta_i - \theta_j))$$
$$= \sum_{j=1}^n (G_{ij}(f_i e_j - e_i f_j) + B_{ij}(e_i e_j + f_i f_j))$$
$$= \sum_{j=1}^n (G_{ij}(W_{i'j} - W_{ij'}) + B_{ij}(W_{ij} + W_{i'j'}))$$

The above expressions indicate that the equality constraints of power injection are linear in terms of $W$. If the cost function is quadratic to $P_g$, then the cost function is quadratic to the elements of $W$. Without the rank 1 constraint $\text{rank}(W) = 1$, this problem is a convex problem and a semi-definite programming (SDP) problem.

## C. Relationship of SOCP and SDP decision variables

In AC OPF with SDP relaxation, $W$ is treated as a decision variable and $W$ should be semi-definite. For SOCP OPF

decision variables ($c_{ij}$ and $s_{ij}$), the following relationship should hold:

$$c_{ij} = e_i e_j + f_i f_j = W_{ij} + W_{i'j'}, \quad (12a)$$
$$s_{ij} = e_i f_j - e_j f_i = W_{ij'} - W_{ji'}, \quad (12b)$$
$$c_{ii} = e_i^2 + f_i^2 = W_{ii} + W_{i'i'}. \quad (12c)$$

For every $c_{ij}$, $s_{ij}$ and $c_{ii}$, we can express them to be the Frobenius product related to the matrix $W$. For example, for a three-bus system with every two buses connected, we have

$$c_{11} = \underbrace{\begin{bmatrix} 1 & 0 & 0 & 0 & 0 & 0 \\ 0 & 0 & 0 & 0 & 0 & 0 \\ 0 & 0 & 0 & 0 & 0 & 0 \\ 0 & 0 & 0 & 1 & 0 & 0 \\ 0 & 0 & 0 & 0 & 0 & 0 \\ 0 & 0 & 0 & 0 & 0 & 0 \end{bmatrix}}_{A_1} \bullet W = \text{Trace}(A_1 W^T) \quad (13)$$

where $\bullet$ denotes Frobenious product. For the test system with 3 buses and 3 lines, we will have 9 variables

$$z = \begin{bmatrix} c_{11} & c_{22} & c_{33} & c_{12} & c_{13} & c_{23} & s_{12} & s_{13} & s_{23} \end{bmatrix}^T$$

and nine $A_i$s.

### D. SDP separation in [3]

Research of [3] indicates that once a vector $z$ is calibrated by SOCP relaxation, this $z$ should be examined. If $z$ is in the set that can be connected with $W$ as shown in (12), then that means $z$ meets the cycle constraint requirement. Otherwise we should create inequalities as cuts. (12) is further expressed as (14) for each cycle in the power network.

$$\mathcal{S} := \left\{ z \in \mathbb{R}^{2|C|} : \exists \tilde{W} \in \mathbb{R}^{2|C| \times 2|C|} \right. \\ \left. \text{s.t.} - z_l + A_l \bullet \tilde{W} = 0 \quad \forall l \in L, \quad \tilde{W} \succeq 0 \right\}, \quad (14)$$

where $C$ is a cycle, and $L$ is the set of lines in the cycles.

The method to create cuts in [3] is shown as follows. For a given $z^*$, the separation problem over $\mathcal{S}$ can be written as follows,

$$v^* := \min_{\alpha, \lambda} \quad -\alpha^T z^* \quad (15a)$$
$$\text{s.t.} \sum_{l \in L} \lambda_l A_l \succeq 0 \quad (15b)$$
$$\alpha + \lambda = 0 \quad (15c)$$
$$-e \leq \alpha \leq e. \quad (15d)$$

If $v^* < 0$, that means the for $z^*$ is infeasible. The cut generated should be $\alpha^T z \leq 0$.

The principle that leads to the above optimization problem was not given in [3]. Here a brief explanation is offered. The sufficient and necessary condition for a matrix $A \in \mathbb{R}^{n \times n}$ is semidefinite is for all $B \in \mathbb{R}^{n \times n}$ and $B \succeq 0$, $A \bullet B \geq 0$.

The proof is given as follows [10]. If $A \succeq 0$ and $B \succeq 0$,

$$A \bullet B = \text{Tr}(AB) = \text{Tr}(A^{1/2} A^{1/2} B^{1/2} B^{1/2}) \\ = \text{Tr}(A^{1/2} B^{1/2} B^{1/2} A^{1/2}) = \|A^{1/2} B^{1/2}\| \geq 0 \quad (16)$$

The converse is proved as follows. If $A \in \mathbb{R}^{n \times n}$ and $A$ is symmetric, and $A \bullet B \geq 0$ for any $B \succeq 0$, then let $x \in \mathbb{R}^n$ and set $B = xx^T$. Then,

$$0 \leq A \bullet B = \text{Tr}(Axx^T) = \sum a_{ij} x_i x_j = x^T Ax. \quad (17)$$

This shows that $A \succeq 0$.

The objective of the SDP separation in [3] is to find a vector $\alpha \in \mathbb{R}^n$ that can make $\alpha^T z_0 \geq 0$ while make $\alpha^T z \leq 0$, where $z \in \mathcal{S}$. The constraint can be further written as

$$\sum_l \alpha_l z_l = (\sum_l \alpha_l A_l) \bullet W \leq 0 \quad (18)$$

The sufficient and necessary condition for the above relation to be true is to have

$$-\sum_l \alpha_l A_l \succeq 0. \quad (19)$$

Therefore, the cut creation can be written in the same format as that in [3].

$$\max_\alpha \quad \alpha^T z_0 \quad (20a)$$
$$\text{st.} - \sum_l \alpha_l A_l \succeq 0 \quad (20b)$$

## III. GENERATING SDP CUTS

### A. Generating SDP Cuts based on Least Square Estimation

The SDP separation method described in [3] assumes that there exists a non-zero vector $\alpha$ to generate a cut expressed in a linear inequality. What if such a vector does not exist? We can easily imagine a two-dimensional case where such cut does not exist. Fig. 1 shows an example that a linear cut cannot separate $Z_0$ from $\mathcal{S}$.

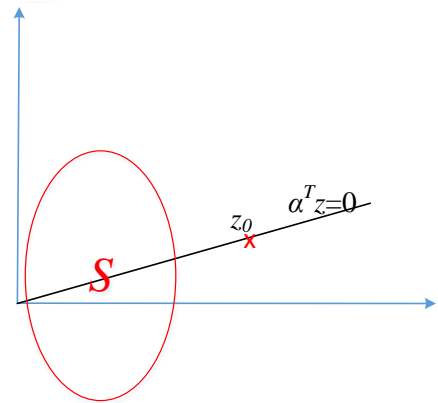

Fig. 1. A case when $\alpha$ that can make $\alpha^T z \leq 0 \quad \forall z \in \mathcal{S}$ and $\alpha^T z_0 > 0$ cannot be found.

In this section, LSE based method is used to find the cuts. The philosophy is explained by Fig. 2. Suppose we have a given $z_0$ and the set we look for is $S$, how do we generate cuts to get rid of $z_0$ and reduce the feasible region and the search space? We use the following method as shown in Fig. 2.

First, we will find the shortest distance from $z_0$ to the set $S$, $z^*$ is the corresponding point found in $S$. The we generate a line that is orthogonal to $z_0 - z^*$. Due to the orthogonality, the vector $z - z_0$ for any $z$ located on the line, is orthogonal to the vector $z_0 - z^*$. Therefore, their inner product is zero.

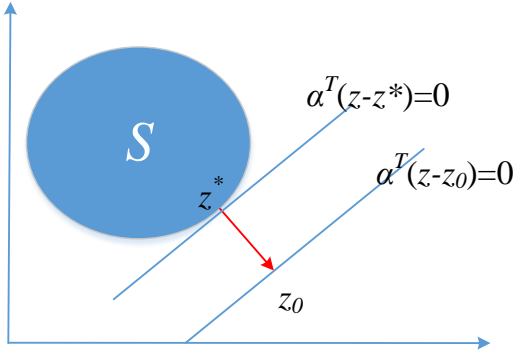

Fig. 2. The vector related to the minimum distance from $z_0$ to $\mathcal{S}$ can be used to create a cut. The cut generated will be $\alpha^T(z - z^*) \leq 0$.

$$(z_0 - z^*)^T(z - z_0) = 0 \tag{21}$$

for any $z$ located on the line. Hence this line is defined as $\alpha^T(z - z_0) = 0$, where $\alpha$ is $(z_0 - z^*)$.

The set $\mathcal{S}$ is now located at the left of the line. Therefore the cut is generated as

$$\alpha^T(z - z_0) \leq 0, \tag{22}$$

where $z$ is the variable.

The task left us is to find $z^*$. This can be done through minimizing the distance from $z_0$ to $z$ where $z \in S$. The formulation is as follows.

$$\min_{z} \quad \|z_0 - z\|^2 \tag{23a}$$
$$st. \quad z_l = \text{Trace}(A_l W^T), \text{ for all } l \in \mathcal{L} \tag{23b}$$
$$W \succeq 0 \tag{23c}$$

where $\mathcal{L}$ is the set of lines. The optimal solution is notated as $z^*$.

If the norm of $z_0 - z^*$ is zero, that means $z_0$ belongs to the SDP set $\mathcal{S}$ and $z_0$ meets the requirement of cycle constraint. The cut generated by (22) is a neutral cut, i.e., when $z = z_0$, $\alpha^T(z - z_0) = 0$. To have a deep cut that can exclude $z_0$, we will use $z^*$, the optimal solution from (23) to generate cut. The cut is expressed as follows.

$$\alpha^T(z - z^*) \leq 0 \tag{24}$$

where $z$ is the variable.

(23) requires to find a $W$ of a large size. This problem will be reduced into small-size optimization problems using cycle basis. That is, for every cycle basis in the system, we will conduct a minimization problem (23) and generate a cut if the vector $z$ generated by SOCP relaxation violates the cycle constraint. Cycle basis identification algorithm in [11] is used to identify the cycle basis.

(23) will be replaced by the following minimization problem for every cycle. For a particular cycle $C$, the edge set is $\mathcal{L}$, the minimization problem is written as follows.

$$\min_{z} \quad \|z_0 - z\|^2 \tag{25a}$$
$$st. \quad z_l = \text{Trace}(A_l \tilde{W}^T), \text{ for all } l \in \mathcal{L} \tag{25b}$$
$$\tilde{W} \succeq 0 \tag{25c}$$

where $z_0$ is part of the vector related to the cycle $C$, $z$ is the decision variable related to the cycle and $\tilde{W} \in \mathbb{R}^{2|C| \times 2|C|}$.

The cut is expressed as follows for each cycle.

$$(z_0 - z^*)^T(z - z^*) \leq 0 \tag{26}$$

### B. Generating SDP Cuts based on feasibility cuts

We further examine how to generate feasibility cuts for a given infeasible decision variable. The basic concept of feasibility cut [12] is first described. Then we examine the feasibility problem presented in (14).

For an inequality constraint $f(x) \leq 0$, where $x \in \mathbb{R}^n$ and $f : \mathbb{R}^n \to \mathbb{R}$ is a convex function, notate the feasible region as $\mathcal{X}$. If $x$ is not feasible and makes $f(x) > 0$, we can find the following relationship based on the convexity of the function $f$:

$$f(z) \geq f(x) + g^T(z - x) \tag{27}$$

where $g$ is a gradient or subgradient. Any feasible $z$ should satisfies the inequality $f(z) \leq 0$. Therefore, we can generate a deep cut as

$$f(x) + g^T(z - x) \leq 0 \tag{28}$$

For the feasibility problem described in (14), we will first come up with an inequality constraint. Note that if $z$ is infeasible, then the following relationship is true.

$$\sum_l (z_l - A_l \bullet \tilde{W})^2 > 0 \tag{29}$$

for any semi-definite $\tilde{W}$.

Since the above condition requires to examine any $\tilde{W}$, the underlying requirement is to examine the following optimization problem:

$$\min_{\tilde{W}} \quad \sum_l (z_l - A_l \bullet \tilde{W})^2 \tag{30a}$$
$$\text{subject to} \quad \tilde{W} \succeq 0 \tag{30b}$$

For a feasible $z$, then

$$\sum_l (z_l - A_l \bullet \tilde{W})^2 \leq 0, \tag{31}$$

for any semidefinite $\tilde{W}$.

The inequality constraint is identified as

$$f(z) = \{\min_{\tilde{W}} \sum_l (z_l - A_l \bullet \tilde{W})^2 \leq 0 | \tilde{W} \succeq 0\}, \tag{32}$$

The above inequality can be further written as

$$f(z) = \min_{\tilde{W}} \quad (z - z^*)^T(z - z^*) \leq 0 \tag{33}$$

where $z_l^* = A_l \bullet \tilde{W}$.

If $z_0$ is infeasible, the feasibility cut should be
$$g^T(z-z_0) + f(z_0) \leq 0 \qquad (34)$$
where $g$ is the gradient. If $g$ is evaluated at $z_0$, then
$$g = 2(z_0 - z^*).$$
In addition, we can find the expression of $f(z_0)$ as
$$f(z_0) = (z_0 - z^*)^T(z_0 - z^*).$$
(34) becomes
$$2(z_0 - z^*)^T(z - z_0) + (z_0 - z^*)^T(z_0 - z^*) \leq 0 \qquad (35)$$
$$\Longrightarrow (z_0 - z^*)^T(z - \frac{z_0 + z^*}{2}) \leq 0. \qquad (36)$$
where
$$z^* = \arg\min\{\min_z \|z_0 - z^*\|^2 \,\big|\, z_l^* = A_l \bullet \tilde{W}, \tilde{W} \succeq 0\}.$$

The feasibility cut is illustrated in the same figure as the LSE-based cut shown in Fig. 3.

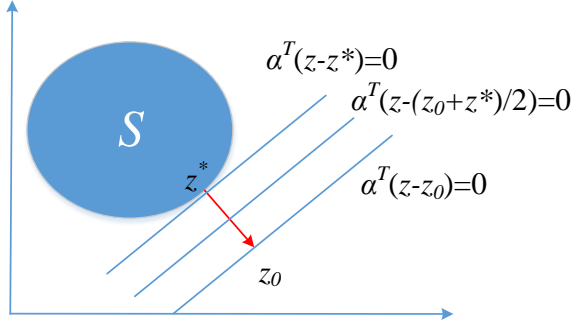

Fig. 3. A neutral cut and two deep cuts. $\alpha^T(z-z^*) \leq 0$ is based on LSE method, $\alpha^T(z - \frac{z_0+z^*}{2}) \leq 0$ is based on feasibility cut. It can be observed that $\alpha^T(z-z^*) \leq 0$ creates the best cut.

### C. Another feasibility cut

We further examine the febrility cut method to arrive at the same cut derived from LSE. Instead of using the square of the norm of $z - z^*$, we will use norm of $z - z^*$. The inequality that decides SDP feasibility is expressed as
$$f(z) = \min_{\tilde{W}} \sqrt{(z-z^*)^T(z-z^*)} \leq 0 \qquad (37)$$
where $z_l^* = A_l \bullet \tilde{W}$.

The gradient of $f(z)$ can be found as:
$$g(z) = \frac{\partial f}{\partial z} = \frac{z - z^*}{\sqrt{(z-z^*)^T(z-z^*)}} \qquad (38)$$

Evaluated at $z_0$, $g(z)$ becomes
$$g(z_0) = \frac{\partial f}{\partial z} = \frac{z_0 - z^*}{\sqrt{(z_0-z^*)^T(z_0-z^*)}} = \frac{z_0 - z^*}{f(z_0)} \qquad (39)$$

Therefore the feasibility cut is
$$g^T(z-z_0) + f(z_0) \leq 0$$
$$\Rightarrow (z_0 - z^*)^T(z - z_0) + f^2(z_0) \leq 0$$
$$\Rightarrow (z_0 - z^*)^T(z - z^*) \leq 0 \qquad (40)$$

Using feasibility cut, we have successfully proved that the LSE-cut is indeed a feasibility cut.

### D. Example SDP cuts for a four-bus system

In this subsection, we use a four-bus system shown in Fig. 4 to illustrate the SDP cuts.

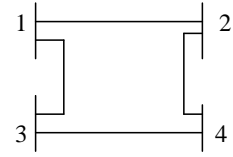

Fig. 4. A four-bus system.

Three sets of decision variables $z$ are generated. The first set is a feasible set that meats the cycle constraint requirement.

1) The first set is generated based on four voltage phasors $V_i \angle \theta_i, i = 1, \cdots 4$.
2) The second set generated by using a set of angle differences ($\theta_{ij}$) that violates the cycle constraint.
3) The third set violates the constraint $c_{ij}^2 + s_{ij}^2 = c_{ii}s_{ii}$.

Table I lists the SOCP decision variables $z_0$, the optimal $z^*$ of the LSE problem, the cut vector $\alpha$. It can be found that for Set 1 and Set 3, the LSE-method will not generate a cut ($\alpha = 0$). For Set 2 (cycle constraint violation), LSE-method can generate a cut ($\alpha \neq 0$).

We also use the three sets of data to generate cuts according to the method in [3]. In all cases, $\alpha = 0$ while $\alpha^T z_0$ is at the order of $10^{-11}$ or approximately zero. SDP separation method [3] cannot generate a cut to exclude an infeasible set (Set 2) that violates cycle constraints.

**Remarks**: The SDP cut generation algorithm has been demonstrated to be effective to identify the infeasible set due to cycle constraint violation. It cannot identify the equality constraint violation related to voltage magnitudes. Our proposed algorithm works better than that in [3]. One reason that the SDP separation in [3] could fail is that the proposed *homogeneous* inequality does not exist. On the other hand, the SDP cut proposed in this paper is *affine* inequality.

### E. Limitation of LSE-SDP Cuts

It has been shown in the previous four-bus example that the proposed SDP cut algorithm cannot detect certain infeasible scenarios (e.g., $c_{ii}$ are disturbed). We will further investigate this scenario to show the generality.

Assume that we have a set of feasible solution $c_{ij}$ and $s_{ij}$ for a cycle. For every $c_{ii}$, we will disturb their values by adding a term $\Delta c_{ii}$. The new set is notated as $c'_{ii}$.
$$c'_{ii} = c_{ii} + \Delta c_{ii}$$
$$= A_i \bullet \tilde{W} + A_i \bullet \Delta W \qquad (41)$$
where $\Delta W$ is a matrix with all elements zero except $i$-th diagonal element ($\Delta W_{ii} = \Delta c_{ii}$).

TABLE I
COMPARISON OF CUT GENERATION FOR TWO SETS OF DECISION VARIABLES (SET 1: FEASIBLE, SET 2: CYCLE CONSTRAINTS VIOLATED, SET 3: CONIC CONSTRAINT VIOLATED)

|  | Set 1 | | | Set 2 | | | Set 3 | | |
| --- | --- | --- | --- | --- | --- | --- | --- | --- | --- |
|  | $z_0$ | $z^*$ | $\alpha$ | $z_0$ | $z^*$ | $\alpha$ | $z_0$ | $z^*$ | $\alpha$ |
| $c_{11}$ | 0.8100 | 0.8100 | -0.00 | 0.8100 | 0.9223 | -0.1123 | 0.9100 | 0.9100 | -0.00 |
| $c_{22}$ | 1.2100 | 1.2100 | -0.00 | 1.2100 | 1.2992 | -0.0892 | 1.3100 | 1.3100 | -0.00 |
| $c_{33}$ | 1.2100 | 1.2100 | -0.00 | 1.2100 | 1.2992 | -0.0892 | 1.3100 | 1.3100 | -0.00 |
| $c_{44}$ | 1.2100 | 1.2100 | -0.00 | 1.2100 | 1.3081 | -0.0981 | 1.2100 | 1.2100 | -0.00 |
| $c_{12}$ | 0.7000 | 0.7000 | 0.00 | 0.9146 | 0.8248 | 0.0898 | 0.7000 | 0.7000 | -0.00 |
| $c_{13}$ | 0.8574 | 0.8574 | 0.00 | 0.8574 | 0.7826 | 0.0747 | 0.8574 | 0.8574 | -0.00 |
| $c_{24}$ | 1.0781 | 1.0781 | 0.00 | 1.1179 | 0.9856 | 0.1323 | 1.0781 | 1.0781 | -0.00 |
| $c_{34}$ | 1.1836 | 1.1835 | 0.00 | 0.8556 | 0.8090 | 0.0466 | 1.1836 | 1.1836 | -0.00 |
| $s_{12}$ | 0.7000 | 0.7000 | 0.00 | 0.3789 | 0.2694 | 0.1094 | 0.7000 | 0.7000 | -0.00 |
| $s_{13}$ | 0.4950 | 0.4950 | 0.00 | -0.4950 | -0.3748 | -0.1202 | 0.4950 | 0.4950 | -0.00 |
| $s_{24}$ | -0.5493 | -0.5493 | -0.00 | -0.4630 | -0.4674 | 0.0043 | -0.5493 | -0.5493 | 0.00 |
| $s_{34}$ | -0.2516 | -0.2516 | -0.00 | -0.8556 | -0.7318 | -0.1238 | -0.2516 | -0.2516 | 0.00 |

Since $A_i$ is a diagonal matrix with only $i$-th and $(i+|C|)$-th diagonal elements 1, we can further let $\Delta W$ be the following diagonal matrix.

$$\Delta W_{11} = \Delta c_{11},$$
$$\vdots$$
$$\Delta W_{ii} = \Delta c_{ii},$$
$$\vdots$$
$$\Delta W_{2|C|,2|C|} = \Delta c_{2|C|,2|C|}$$

Therefore, for any $c_{jj}$, the disturbed $c'_{jj}$ will all share the same $\Delta W$.

$$\begin{aligned} c'_{jj} &= c_{jj} + \Delta c_{jj} \\ &= A_j \bullet \tilde{W} + A_j \bullet \Delta W \\ &= A_j \bullet \tilde{W}' \end{aligned} \quad (42)$$

The above analysis shows that if $c_{ij}^2 + s_{ij}^2 < c_{ii}c_{jj}$, the SDP cut algorithm will not be able to detect the infeasibility. Instead, SOCP AC OPF takes the task to enforce $c_{ij}^2 + s_{ij}^2 = c_{ii}c_{jj}$ by the objective function minimization procedure. Minimizing the total active power loss or minimizing the cost of generators all help to enforce $c_{ij}^2 + s_{ij}^2 = c_{ii}c_{jj}$ under the conditions, *e.g.*, voltage upper bounds not binding [8].

## IV. CASE STUDIES

Experimental Setting: All computations are conducted in MATLAB. MATPOWER [13] is used to find upper bound while CVX toolbox [14] is used to carry out convex programming problem solving. Cycle basis is identified using the algorithm in [11].

### A. Three-bus test system

This test case comes from NESTA v0.4.0 archive [15]. The system is a three-bus system consisting of one cycle. After adding five SDP cuts, the solution from SOCP AC OPF is in the SDP set. Table II lists the gap before adding SDP cuts and after adding SDP cuts. Value of the minimization problem that generates SDP cuts is plotted for five iterations.

TABLE II
PERCENTAGE GAP OF A 3-BUS CASE

| Test case | MATPOWER ($/h) | SOCP | SDP cuts |
| --- | --- | --- | --- |
| nesta_case3_lmbd | 5812.6 | 1.67 | 1.27 |

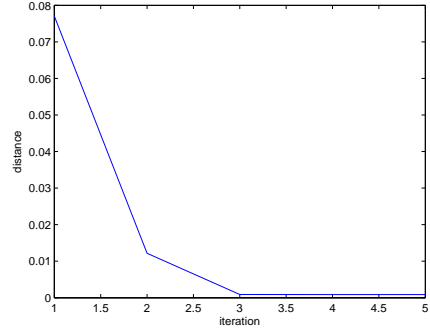

Fig. 5. Three-bus test case: $z_0$ to $\mathcal{S}$ distance over iteration.

### B. Five-bus PJM system

This test case also comes from NESTA archive. The five-bus system has two cycles. Therefore, each time, two cuts will be generated. After five iterations, the gap reduces from 14.5% to 8.8%.

TABLE III
PERCENTAGE GAP OF A 5-BUS CASE

| Test case | MATPOWER ($/h) | SOCP | SDP cuts |
| --- | --- | --- | --- |
| nesta_case5_pjm | 17552 | 14.54 | 8.88 |

### C. 30-bus test case

Nesta's IEEE 30-bus system has been reported in [3], [4] to have a large gap for SOCP relaxations and quadratic relaxations. This system has 12 cycles. SOCP relaxation gives a set of decision variables. With the set, 12 minimum distance problems are solved and 12 SDP cuts are generated. The SCOP problem with 12 cuts is solved again and the gap is reduced significantly to 1.71%. After five iterations, the gap is reduced to 0.42%. Fig. 7 gives the 12 distances over iteration.

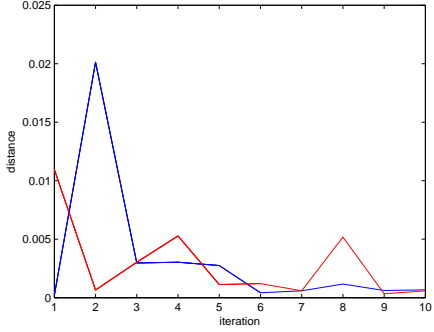

Fig. 6. Five-bus test case: $z_0$ to $\mathcal{S}$ distance over iteration.

TABLE IV
PERCENTAGE GAP OF A 30-BUS CASE

| Test case | MATPOWER ($/h) | SOCP | SDP cuts |
|---|---|---|---|
| nesta_case30_ieee | 204.76 | 15.88 | 0.42 |

*D. 118-bus system and 300-bus system*

Two systems with hundreds of buses are tested. These two systems are from MATPOWER's database. For each case, two iterations are conducted. The results are presented in Table V. In both cases, the gap is reduced significantly.

TABLE V
PERCENTAGE GAP OF 118-BUS AND 300-BUS SYSTEMS

| Test case | MATPOWER ($/h) | cycles | SOCP | SDP cuts |
|---|---|---|---|---|
| 118-bus | 130395 | 62 | 0.25 | 0.1 |
| 300-bus | 719740 | 110 | 0.15 | 0.06 |

A number of Nesta cases have been tested and the results are listed in Table VI lists. Gaps are compared for SOCP relaxation, SOCP relaxation with LSE-based SDP cuts, and quadratic programming (QC) relaxation [4]. The gaps for SOCP relaxation, SOCP relaxation with LSE-based SDP cuts are computed from this project. The gaps of QC relaxation come from [4].

## V. CONCLUSION

This paper presents LSE-based approach to find affine inequalities for AC OPF with SOCP relaxation. The affine inequalities serve as SDP cuts to reduce the feasible region and get rid of SOCP solutions outside of the feasible region of SDP relaxation. Least square estimation-based SDP cuts have been demonstrated to be effective to cut infeasible solutions for AC OPF with SOCP relaxations. This method has been tested on variety of cases to demonstrate its effectiveness.

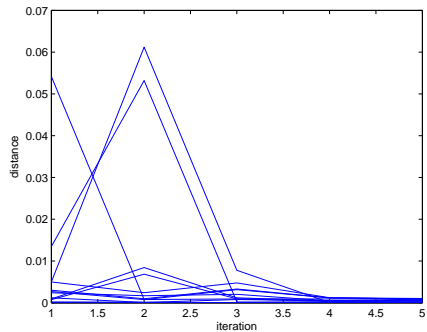

Fig. 7. 30-bus test case: $z_0$ to $\mathcal{S}$ distance over iteration.

TABLE VI
OPF RESULTS FOR ALL TEST CASES.

| Test case | SOCP | SDP cuts | QC [4] |
|---|---|---|---|
| nesta_case3_lmbd | 1.67 | 1.27 | 1.24 |
| nesta_case4_gs | 0.00 | 0.00 | n.a. |
| nesta_case5_pjm | 14.54 | 8.88 | 14.54 |
| nesta_case6_ww | 0.63 | 0.02 | n.a. |
| nesta_case9_wscc | 0.00 | 0.00 | n.a. |
| nesta_case14_ieee | 0.11 | 0.00 | n.a. |
| nesta_case29_edin | 0.14 | 0.05 | n.a. |
| nesta_case30_as | 0.06 | 0.00 | n.a. |
| nesta_case30_fsr | 0.39 | 0.19 | n.a. |
| nesta_case30_ieee | 15.88 | 0.42 | 15.64 |
| nesta_case39_epri | 0.05 | 0.02 | n.a. |
| nesta_case57_ieee | 0.06 | 0.01 | n.a. |
| nesta_case118_ieee | 2.22 | 1.58 | 1.72 |
| nesta_case162_ieee | 2.07 | 1.82 | 4.00 |
| nesta_case300_ieee | 1.33 | 0.79 | 1.17 |

**Zhixin Miao** (S'00 M'03 SM'09) received the B.S.E.E. degree from the Huazhong University of Science and Technology,Wuhan, China, in 1992, the M.S.E.E. degree from the Graduate School, Nanjing Automation Research Institute, Nanjing, China, in 1997, and the Ph.D. degree in electrical engineering from West Virginia University, Morgantown, in 2002.

Currently, he is with the University of South Florida (USF), Tampa. Prior to joining USF in 2009, he was with the Transmission Asset Management Department with Midwest ISO, St. Paul, MN, from 2002 to 2009. His research interests include power system stability, microgrid, and renewable energy.

**Lingling Fan** received the B.S. and M.S. degrees in electrical engineering from Southeast University, Nanjing, China, in 1994 and 1997, respectively, and the Ph.D. degree in electrical engineering from West Virginia University, Morgantown, in 2001. Currently, she is an Associate Professor with the University of South Florida, Tampa, where she has been since 2009. She was a Senior Engineer in the Transmission Asset Management Department, Midwest ISO, St. Paul, MN, form 2001 to 2007, and an Assistant Professor with North Dakota State University, Fargo, from 2007 to 2009. Her research interests include power systems and power electronics. Dr. Fan serves as a technical program committee chair for IEEE Power System Dynamic Performance Committee and an editor for IEEE Trans. Sustainable Energy.

**Hossein Ghassempour Aghamolki** received B. S. degree in Electrical Engineering from Power & Water University of Technology, Tehran, Iran, in 2007 and the M. S. in Electrical Engineering from University of Mazandaran, Babol, Iran, in 2010. Currently, He is a Ph.D. student at University of South Florida (USF), Tampa FL. He started his Ph.D. study in Fall 2013. His research interests include state estimation and system identification, dynamic modelling, optimization and power systems.

**Bo Zeng (M'11)** received the Ph.D. degree in industrial engineering from Purdue University, West Lafayette, IN, USA. He currently is an Assistant Professor with the Department of Industrial Engineering, University of Pittsburgh, Pittsburgh, PA, USA. His research interests include polyhedral study and algorithms for stochastic and robust mixed integer programs, coupled with applications in power and logistics systems. He is a member of IIE and INFORMS.